# ABOUT BERNOULLI'S NUMBERS

Prof. Mihàly Bencze
Department of Mathematics
Áprily Lajos College
Braşov, Romania

Florentin Smarandache, Ph D
Full Professor
Chair of Department of Math & Sciences
University of New Mexico
200 College Road
Gallup, NM 87301, USA
E-mail: smarand@unm.edu

**Abstract**.
Many methods to compute the sum of the first $n$ natural numbers of the same powers (see [4]) are well known.
In this article we present a simple proof of the method from [3].

**Introduction**.
The Bernoulli's numbers are defined by

(1) $\quad B_n = \dfrac{-1}{n+1}\left(C_{n+1}^0 B_0 + C_{n+1}^1 B_1 + \ldots + C_{n+1}^{n-1} B_{n-1}\right)$

where $B_0 = 1$. It is known that $B_{n+1} = 0$ if $n \geq 1$. By calculation we find that:

(2) $\quad B_1 = -\dfrac{1}{2},\ B_2 = \dfrac{1}{6},\ B_4 = -\dfrac{1}{30},\ B_6 = \dfrac{1}{42},\ B_8 = -\dfrac{1}{30},\ B_{10} = \dfrac{5}{66},$

$\quad B_{12} = -\dfrac{691}{2730},\ B_{14} = \dfrac{7}{6},\ B_{16} = -\dfrac{3617}{510},\ B_{18} = \dfrac{43867}{798},\ B_{20} = -\dfrac{174611}{330},$

$\quad B_{22} = \dfrac{854513}{138},\ B_{24} = -\dfrac{236364091}{2730}$, etc.

Let $S_n^k = 1^k + 2^k + \ldots + n^k$ the sum of the first $n$ natural numbers which have the same power.

**Theorem.**

(3) $\quad S_n^k = \dfrac{1}{k+1}\left(n^{k+1} + \dfrac{1}{2}C_{k+1}^1 n^k + C_{k+1}^2 B_2 n^{k-1} + \ldots + C_{k+1}^k B_n n\right)$

*Proof:* (1) can be written as:

(4) $$\sum_{i=0}^{n} C_{n+1}^{i} B_{i} = 0, \ n \geq 1.$$

If $P(x) = \sum_{i=0}^{k} C_{k+1}^{i} B_i x^{k+1-i}$,

then

$$P(n+1) - P(n) = \sum_{i=0}^{k} C_{k+1}^{i} B_i \left( (n+1)^{k+1-i} - n^{k+1-i} \right) = \sum_{i=0}^{k} C_{k+1}^{i} B_i \left( \sum_{j=1}^{k+1-i} C_{k+1-i}^{j} n^{k+1-i-j} \right).$$

Let $A_t$ be the coefficients of $n^{k-t}$, where $t \in \{0, 1, ..., k\}$.

$$A_t = \sum_{i=0}^{t} C_{k+1}^{i} C_{k+1-i}^{t+1} B_i = C_{k+1}^{t+1} \left( \sum_{i=0}^{t} C_{t+1}^{i} B_i \right).$$

If $n \geq 1$, then $A_t = 0$, only $A_0 = C_{k+1}^{1}$.

Because of these $P(n+1) - P(n) = C_{k+1}^{1} n^k$. Using this

$$\sum_{i=0}^{n-1} i^k = \frac{1}{k+1} \sum_{i=0}^{n-1} (P(i+1) - P(i)) = \frac{1}{k+1} P(n),$$

because $P(0) = 0$. Then $S_n^k = \frac{1}{k+1} P(n) + n^k$. From here one obtains (3).

**Note**. From the previous result we can also find the formula

$$S_n^k = \frac{1}{k+1} P(n+1).$$

Using this, we find the following equalities:

$$S_n^0 = n, \ S_n^1 = \frac{1}{2} n(n+1), \ S_n^2 = \frac{1}{6} n(n+1)(2n+1), \ S_n^3 = \frac{1}{4} n^2 (n+1)^2,$$

$$S_n^4 = \frac{1}{30} n(n+1)(2n+1)(3n^2 + 3n - 1), \ S_n^5 = \frac{1}{12} n^2 (n+1)^2 (2n^2 + 2n - 1),$$

$$S_n^6 = \frac{1}{42} n(n+1)(2n+1)(3n^4 + 6n^3 - 3n + 1),$$

$$S_n^7 = \frac{1}{24} n^2 (n+1)^2 (3n^4 + 6n^3 - n^2 - 4n + 2),$$

$$S_n^8 = \frac{1}{90} n(n+1)(2n+1)(5n^6 + 15n^5 + 5n^4 - 15n^3 - n^2 + 9n - 3),$$

$$S_n^9 = \frac{1}{20} (2n^{10} + 10n^9 + 15n^8 - 14n^6 + 10n^4 - 3n^2),$$

$$S_n^{10} = \frac{1}{66} (6n^{11} + 33n^{10} + 55n^9 - 66n^7 + 66n^5 - 33n^3 + 5n),$$

$$S_n^{11} = \frac{1}{24} (2n^{12} + 12n^{11} + 22n^{10} - 33n^8 + 44n^6 - 33n^4 + 10n^2),$$

$$S_n^{12} = \frac{1}{2730} (210n^{13} + 1365n^{12} + 3630n^{11} - 4935n^9 + 115n^8 +$$

$$+ 9640n^7 + 1960n^6 - 5899n^5 + 35n^4 + 4550n^3 + 1382n^2 - 691n), \text{ etc.}$$

**Problems:**
1) Using the mathematical induction on the base of (1), we prove that
   $B_{2n+1} = 0$, if $n \geq 1$.
2) Prove that $S_n^k$ is divisible by $n(n+1)$.
3) Prove that $S_n^{2k+1}$ is divisible by $n^2(n+1)^2$.
4) Determine those natural numbers $n, k$ for which $S_n^{2k}$ is divisible $n(n+1)(2n+1)$.
5) Detach in parts the sums $S_n^9$, $S_n^{10}$, $S_n^{11}$, $S_n^{12}$.
6) Using (2), (3), compute the sums $S_n^{13},...,S_n^{21}$.